\documentclass{article}

\usepackage[a4paper, left=3.3cm,top=3.3cm,right=3.3cm,bottom=3.7cm]{geometry}
\usepackage{concmath}
\usepackage[T1]{fontenc}

\usepackage{graphicx}
\usepackage{mathtools}
\usepackage{psfrag}
\usepackage{subfig}

\usepackage{xspace}
\usepackage{pgfplots}
\usepackage{siunitx}

\usepackage{enumitem}
\setitemize{label=\scriptsize{$\blacksquare$},topsep=0em}
\setenumerate{label=(\alph*), topsep=0em}

\newcommand\inner[2]{\left\langle  #1,#2 \right\rangle}
\newcommand{\kl}[1]{\left(#1\right)}

\newcommand{\trans}{\mathsf{T}}

\newcommand{\R}{\mathbb R}

\newcommand{\N}{\mathbb N}

\def\plus{{\boldsymbol{\texttt{+}}}}

\newcommand{\al}{\alpha}

\newcommand{\ran}{\operatorname{ran}}

\newcommand{\rr}{\mathbf{r}}
\newcommand{\rs}{\mathbf{s}}
\newcommand{\cx}{x}
\newcommand{\cy}{y}
\newcommand{\noise}{\xi}

\newcommand\norm[1]{{\left\Vert#1\right\Vert}}

\newcommand{\Ao}{\mathbf  A}                         
\newcommand{\Bo}{\mathbf  B}                         
\newcommand{\Ro}{\mathbf R}    
\newcommand{\Pro}{\mathbf P}    %
\newcommand{\Uo}{\mathbf U}
\newcommand{\Sigmao}{\mathbf{\Sigma}}
\newcommand{\Vo}{\mathbf{V}}

\newcommand{\Id}{\operatorname{Id}}                         

\newcommand{\nun}{{\mathbf \Phi}}  
\newcommand{\rmd}{\mathrm d}

\usepackage{soul}
\usepackage{xcolor}
\colorlet{lred}{red!40}
\colorlet{lgreen}{green!40}
\colorlet{lblue}{blue!40}

\definecolor{lime}{RGB}{34,139,34}
\definecolor{goyel}{RGB}{90,215,90}
\definecolor{blor}{RGB}{250,188,81}
\definecolor{blgr}{RGB}{210,105,30}
\definecolor{orred}{RGB}{255,165,0}
\definecolor{orbl}{RGB}{90,215,90}
\definecolor{dot}{RGB}{0,0,95}

\def\mirrord at (#1,#2){\filldraw[color=orred](#1-0.07,#2)--(#1-0.07,#2-0.2)--(#1-0.12,#2-0.2)--(#1,#2-0.35)--(#1+0.12,#2-0.2)--(#1+0.07,#2-0.2)--(#1+0.07,#2)--cycle;}

\newcommand{\equic}[3][1 cm]{
  \foreach \i in {1,...,#2} {
    \coordinate (N\i) at ({#1*cos(\i*180/(#2+1))+#3cm},{-#1*sin(\i*180/(#2+1))});
    \fill[black] (N\i) circle (0.1 cm);
}
}

\usepackage{amsmath,amsfonts,amssymb}
\usepackage{graphicx}
\usepackage[colorlinks=true, allcolors=blue]{hyperref}

\usepackage{authblk}
\title{Deep Learning of truncated singular values for limited view photoacoustic tomography}

\author{Johannes Schwab}
\author{Stephan Antholzer}

\affil{Department of Mathematics, University of Innsbruck\authorcr
Technikerstrasse 13, 6020 Innsbruck, Austria\authorcr
{\tt \{johannes.schwab, stephan.antholzer\}@uibk.ac.at}}

\author{Robert Nuster}
\author{G\"unther Paltauf}

\affil{Department of Physics, Universit\"at Graz\authorcr
Universitaetsplatz 5, Graz, Austria\authorcr
E-mail: {\tt \{ro.nuster, guenther.paltauf\}@uni-graz.at }}

\author{Markus Haltmeier}

\affil{Department of Mathematics, University of Innsbruck\authorcr
Technikerstrasse 13, 6020 Innsbruck, Austria\authorcr
E-mail: {\tt markus.haltmeier@uibk.ac.at}}

\date{Januar 19, 2019}

\begin{document}
\maketitle

\begin{abstract}
We develop a  data-driven regularization method
for the severely  ill-posed problem of photoacoustic image reconstruction from limited view data. Our approach is based on the  regularizing networks  that have been  recently  introduced and analyzed in [J. Schwab, S. Antholzer, and M. Haltmeier. \emph{Big in Japan: Regularizing networks for solving inverse problems (2018)}, arXiv:1812.00965] and consists of two steps.
In the first step, an intermediate reconstruction is performed  by applying  truncated singular value decomposition (SVD).
In order to prevent noise amplification,  only  coefficients   corresponding  to sufficiently large singular values are  used, whereas the remaining  coefficients are set zero. In a second step, a trained deep neural network is  applied to recover the truncated SVD coefficients. Numerical results are presented demonstrating  that the proposed data driven estimation of the truncated singular values significantly improves the pure truncated SVD reconstruction. We point out that proposed reconstruction  framework
can straightforwardly  be applied  to other inverse  problems, where the SVD is either known analytically or can be computed numerically.


\end{abstract}

\section{Introduction}
\label{sec:intro}  

Recently, a considerable amount of research has been done on the use of deep neural networks for various  image reconstruction tasks. Deep convolutional neural networks (CNNs) have shown excellent results  for  deconvolution, inpainting or denoising problems. They have also been successfully used in tomographic reconstruction problems, like full waveform inversion, x-ray tomography, magnetic resonance imaging and photoacoustic tomography (see, for example \cite{adler2017solving,antholzer2018deep,bubba2018learning,han2018framing,gupta2018cnn,kobler2017variational,li2018nett,ye2018deep}).
There exist several different deep learning  approaches to solve inverse problems. This includes CNNs in iterative schemes \cite{gupta2018cnn,kobler2017variational,adler2017solving}, or using a single CNN to improve an initial reconstruction \cite{han2018framing,jin2017deep,schwab2018deep}.
Very recently, we have proven that sequences of suitable neural networks in combination with  classical regularizations lead to  convergent regularization methods \cite{schwab2018big}.
In the present paper, we apply the concept of  regularizing networks to the limited view problem of photoacoustic tomography (PAT).  In particular, we use truncated SVD  as initial regularization method   for  approximating a low frequency  approximation of the photoacoustic  (PA) source. In the  second step,  a deep network is trained  to recover the missing high frequency part.

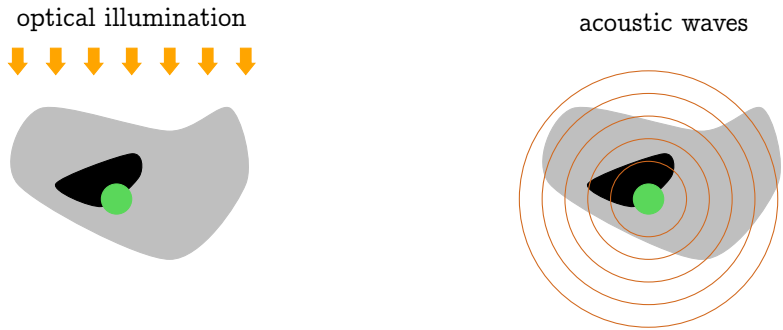
\begin{figure}[htb!]
\centering
\begin{tikzpicture}
\mirrord at (0,1.8);
\mirrord at (0.5,1.8);
\mirrord at (1,1.8);
\mirrord at (1.5,1.8);
\mirrord at (2,1.8);
\mirrord at (2.5,1.8);
\mirrord at (3,1.8);
\draw (1.5,1.9) node[above]{optical illumination};
\filldraw[gray!50] plot[smooth cycle] coordinates {(0,0)(2,-1)(3,0)(2.8,1)(2,0.7)(0.3,1)};
\filldraw[black] plot[smooth cycle] coordinates {(0.5,0)(1.1,-0.3)(1.6,0.1)(1.5,0.4)};
\filldraw[color=orbl] (1.3,-0.2) circle(0.2);
\filldraw[gray!50] plot[smooth cycle] coordinates {(7,0)(9,-1)(10,0)(9.8,1)(9,0.7)(7.3,1)};
\filldraw[black] plot[smooth cycle] coordinates {(7.5,0)(8.1,-0.3)(8.6,0.1)(8.5,0.4)};
\filldraw[color=orbl] (8.3,-0.2) circle(0.2);
\draw[color=blgr] (8.3,-0.2) circle(0.5);
\draw[color=blgr] (8.3,-0.2) circle(0.8);
\draw[color=blgr] (8.3,-0.2) circle(1.1);
\draw[color=blgr] (8.3,-0.2) circle(1.4);
\draw[color=blgr] (8.3,-0.2) circle(1.7);
\draw (4,-2.6) node{};
\draw (8.5,1.9) node[above]{acoustic waves};
\end{tikzpicture}
\caption{\textbf{Basic principles of PAT.}
The absorption of short optical pulses inside a semitransparent sample causes thermoelastic expansion, which in turn induces acoustic pressure waves. The acoustic waves propagate outwards, where they are
measured and used to reconstruct an image of the interior.}\label{fig:pateffect}
\end{figure}

PAT is a tomographic imaging  method based on the generation of acoustic waves induced by pulsed optical  illumination  (see Figure~\ref{fig:pateffect}).   The induced pressure waves are measured outside of the investigated sample  and used for  reconstructing the  PA source. In many practical applications, the pressure data cannot be measured on a surface fully surrounding the sample, which makes the reconstruction problem severely ill-posed.  This implies that exact  inversion methods
are unstable and significantly  amplify noise. Therefore, one has to apply regularization techniques which replace the exact inverse (or pseudo-inverse) by approximate  but stable  inversion methods. In this paper, we  work with a fully discretized model based approach, where  we recover the  expansion coefficients  of the PA source with respect to certain
 basis  functions \cite{dean2012accurate,paltauf2007experimental,paltauf2002iterative,rosenthal2013acoustic,zhang2009effects,schwab2018galerkin}.
The discretized limited view scenario  yields  to the solution of a
severely ill-conditioned system of linear equations.  The ill-conditioning is reflected by the  fact that the singular values
decay rapidly and a large fraction of them is close to zero.

In order to solve the discrete ill-conditioned limited view problem of PAT, we propose the following deep-learning based reconstruction procedure. In a first step, we apply  truncated SVD  to obtain an intermediate  reconstruction formed by the singular components corresponding to sufficiently large singular values. By truncating the small singular values we prevent the reconstruction to amplify noise since small singular values of the forward matrix correspond to large  singular values of the inverse matrix that  amplify noise. For  recovering  the truncated singular components, in a second step, we train a deep CNN that maps the intermediate reconstruction to the truncated part. The actual reconstruction  is the sum of the  intermediate reconstruction and the residual image found by the network. Opposed to strict null-space learning proposed in \cite{schwab2018deep}, not only the part in the kernel of the operator (corresponding to singular value zero) is learned, but also parts corresponding to small singular values. Consequently, the  approach proposed here can be seen as approximate null-space learning.
 Our numerical results demonstrate, that the proposed  deep learning based extension of truncated singular components  significantly improves pure truncated SVD.

\section{Photoacoustic tomography}

\subsection{Continuous model}

Throughout this paper we work with homogeneous medium and  describe the  acoustic wave propagation in PAT by the initial value problem
\begin{align}
\partial_t^2 p(\rr,t)-c^2\Delta_\rr p(\rr,t)&=0\,,\nonumber\\
p(\rr,0)&=f(\rr)\,,\label{eq:waveeq}\\
\partial_t p(\rr,0)&=0 \,.\nonumber
\end{align}
Here $d=2$ or $d=3$ is the spatial dimension, $f \colon \R^d \to \R$
 the PA source, $\rr \in\R^d$  the
 spatial variable, $t\in\R$ the  time variable, $\Delta_\rr$  the spatial Laplacian, and $c$ is the constant speed of sound. We assume that the PA exists inside the region $R \subset \R^d$.
 After temporal rescaling, in the following we assume  normalized speed
 of sound $c=1$.  The inverse problem of PAT consists in reconstructing  $f$ from values of $p(\rs,t)$ for points   $\rs  \in S $ on a measurement surface $S \subset \partial R$ and times $t \geq 0$.
For point-like detectors, the spatial dimension $d=3$ is relevant. In this work, we consider  the case $d=2$, which is relevant for PAT using integrating line detectors \cite{burgholzer2007temporal, paltauf2007photoacoustic}. Using integrating line detectors leads to a set of 2D problems (\eqref{eq:waveeq} for $d=2$), where the initial PA source   corresponds to the projection of the 3D source in a certain direction. By varying  the projection direction  and applying  the inverse Radon transform to the  projection images one can reconstruct the 3D source, although the projection images are already of diagnostic value.

We denote by  $\mathcal{A}$ the operator, that maps a PA source $f$ to the solution $\mathcal{A} f \colon \partial R \times [0,\infty)\rightarrow \R$ restricted to $\partial R$.
In the full data situation, the goal is  to reconstruct  $f$ from noisy versions of $\mathcal{A} f$.
This problem can be solved stably by means  of explicit inversion formulas
known for certain geometries \cite{burgholzer2007temporal,Kun07,finch2007inversion,XW05,finch2004determining,haltmeier2014universal,palamodov2014time,nguyen2014reconstruction}. In practical applications, measurement data are only available  on a surface $S \subsetneq \partial R$ that does not fully enclose the investigated sample. In this case, the inverse problem of PAT becomes severely
ill-posed. \cite{xu2004reconstructions,frikel2015artifacts,barannyk2015artifacts,stefanov2009thermoacoustic}
In this work we consider the limited view problem in  a discrete
setting as  described in the following subsection.

\subsection{Discretization}

For the discretization we assume  that the initial pressure has the following form
\begin{align}\label{eq:approx}
f(\rr) = \sum_{i=1}^{N^2} \cx_i \, \varphi_i(\rr)  \quad \text{ for expansion coefficients $\cx_i \in \R$}  \,.
\end{align}
Here $\varphi_i (\rr) = \varphi(\rr-\rr_i)$ are  translated version of a fixed  basis function $\varphi\colon \R^2\rightarrow \R$, where the centers $\rr_i$ are
arranged on a Cartesian grid.
In particular, we consider so-called Kaiser-Bessel functions which are radially symmetric functions of compact support that are  very popular for tomographic inverse problems \cite{wang2014discrete,wang2012investigation,schwab2018galerkin}. The generalized Kaiser-Bessel functions depend on three parameters and  are defined by
\begin{equation*}
\varphi(\rr) \triangleq
\begin{cases}
\left(\sqrt{1-\norm{\rr}^2/a^2}\right)^m\frac{I_m\left(\gamma \sqrt{1-\norm{\rr}^2/a^2}\right)}{I_m(\gamma)} \quad &\text{ for } \norm{\rr}\geq a\\
0 \quad &\text{ else} \,.
\end{cases}
\end{equation*}
Here $I_m$ is the modified Bessel function of the first kind of order $m\in\N$, $\gamma>0$ the window taper and $a$ the radius of the support.
Since the forward operator $\mathcal{A}$ is linear we have
$ \mathcal{A} f = \sum_{i=1}^{N^2} \cx_i \mathcal{A} \varphi_i$.

In applications, the number of spatial and  temporal measurements   is finite. We  assume that measurements are made with
 sampling points $(\rs_n)_{n=1}^{N_s}$ on  the detection surface $S$
 and temporal sampling points  $(t_j)_{j=1}^{N_t}$ in $[0, T]$.
In our experiments presented below, we consider equidistant spatial
sampling points $\rs_n = (\cos(n\pi/N_s),\sin(n \pi/N_s))$
for $n=1,\ldots, N_s$ and equidistant temporal  sampling points $t_j = jT/N_t$  for  $j=1,\ldots, N_t$.
For  each  basis function $\varphi_i$, we  write the corresponding data
$\mathcal{A} \varphi_i$ evaluated at the discretization points
into the columns of a  matrix, $\Ao_{N_t (n-1)+j, i }  \triangleq   \mathcal{A}\varphi_i(\rs_n,t_j)$.
This results in a  linear system of equations for the expansions coefficients of the PA source:
\begin{equation}\label{eq:linsys}
\text{Recover $\cx \in \R^{N^2}$  from    data } \quad
	\cy = \Ao \cx + \noise \,.
\end{equation}
Here $\Ao \in \R^{N_t N_s\times N^2}$ is the system matrix,
$\noise$   models the error in  the data,
and $\cy \in \R^{N_t N_s}$ is vector of the noisy data  $\cy_{N_t (n-1)+j} = \mathcal{A}f(\rs_n,t_j) + \noise_{N_t (n-1)+j, i} $
given on the measurement points.
Equation \eqref{eq:linsys} is ill-conditioned, reflecting the ill-posedness of the underlying continuous limited view problem.

In the following, we consider the fully discrete problem of recovering the vector $\cx$ of expansion coefficients in \eqref{eq:linsys}.
The PA source can then be recovered by evaluating the expansion \eqref{eq:approx}. Recall that $\varphi_i$ are  translated versions of a fixed  basis function with centers on a Cartesian grid.
We therefore  can arrange the entries of a vector $\cx \in \R^{N^2}$ as
$N \times N$ image in a natural manner.
This image representation will be used when visualizing elements
in $\R^{N^2}$ and for applying the CNNs.

\section{Deep learning of singular value expansion}
\label{sec:reg}

In this section we describe the proposed deep learning method
for solving the discrete linear system \eqref{eq:linsys}.
As it is shown in \cite{schwab2018big}
the combination of a classical regularization method with deep CNNs  yields a regularization method together with quantitative error estimates. Here we  use the truncated SVD decomposition as classical
regularization method and combine it with a deep network
that recovers the truncated coefficients.

\subsection{Truncated SVD reconstruction}

In the following, let $\Ao = \Uo \Sigmao \Vo^\trans$  denote  the singular value decomposition of the forward matrix $\Ao$. Here $\Uo\in \R^{N_t N_s\times N_t N_s}$ and $\Vo\in \R^{N^2 \times N^2}$ are orthogonal matrices, and $\Sigmao \in \R^{N_t N_s\times N}$ is the  diagonal matrix. The tuples of columns $(u_i)_{i=1}^{N_t N_s}$ and $(v_i)_{i=1}^{N^2}$ are orthonormal systems in $\R^{N_t N_s}$ and $\R^{N^2}$, respectively. The diagonal entries $\sigma_1 \geq \sigma_2 \geq \cdots \geq \sigma_{ \min\{N_t N_s, N^2\}} \geq 0$ of $\Sigmao$ are the singular  values of $\Ao$. The number $r \leq \min\{N_t N_s, N^2\}$ of non-vanishing singular values equals the rank of the matrix $\Ao$.
Using the SVD, we have  the explicit formulas for the forward matrix $\Ao$ and its   Moore-Penrose pseudo-inverse $\Ao^\plus $ (also called Moore-Penrose inverse),
\begin{align}\label{eq:svd1}
\Ao \cx &= \sum_{i=1}^r \sigma_i \inner{u_i}{\cx} v_i
\quad  \text{ for }  \cx \in \R^{N^2}
\\ \label{eq:svd2}
\Ao^\plus \cy= \Vo \Sigmao^\plus \Uo^\trans \cy
& = \sum_{i=1}^r \frac{1}{\sigma_i}\inner{\cy}{v_i} u_i
\quad  \text{ for }  \cy \in \R^{N_t N_s}\,.
\end{align}
In particular, for exact data $\cy = \Ao \cx$, \eqref{eq:svd2} gives an explicit solution for
the discrete limited view PAT problem \eqref{eq:linsys}.
In the case of noisy data $\cy = \Ao \cx+\noise$, the small singular values $\sigma_i$ cause severe amplification of the noise coefficients $\inner{\noise}{v_i}$ and typically yields to useless reconstruction results.
To stabilize the inversion problem one has to apply regularization
methods which replace  the exact pseudo-inverse
by a stable approximate inverse.

In this paper, we work with truncated SVD which is a well established regularization method.
Truncated SVD consists of a family $(\Bo_\al)_{\al >0}$ of
approximations of  the pseudo-inverse  defined by
\begin{equation} \label{eq:tsvd}
\Bo_\al(\cy) = \sum_{\sigma_i^2\geq \al} \frac{1}{\sigma_i}\inner{\cy}{v_i} u_i \quad \text{for } \cy\in\R^{N_t N_s}.
\end{equation}
Applying truncated SVD  with  regularization parameter
$\alpha >0$ (which is chosen depending on the noise level) prevents amplification of measurement errors.
In particular, when applied to noisy data  satisfying $\norm{\Ao \cx - \cy } \leq \delta$ for some noise level $\delta>0$,
it satisfies the stability estimate $\norm{\Bo_\al  \cy - \Ao^\plus \Ao  \cx} \leq \delta / \sqrt{\alpha}
+  \norm{ (\Ao^\plus   - \Bo_\al ) \Ao  \cx}$. This implies convergence of truncated SVD
in the  case that the regularization parameter $\al = \alpha(\delta)$ is chosen such that
$\al \to 0$ and $\delta^2/\al \to 0$ as $\delta \to 0$.   Quantitative error estimates require
bounding the approximation term  $\norm{ (\Ao^\plus   - \Bo_\al ) \Ao  \cx}$. Noting that  $(\Ao^\plus   - \Bo_\al ) \Ao  \cx
= \sum_{\sigma_i^2 < \al}  \inner{\cx}{u_i} u_i$ this requires  $\cx$
being  smooth in the sense that the contribution of coefficients corresponding to small singular values is small.

\subsection{Deep Learning of truncated singular values}

According to \eqref{eq:tsvd},    the truncated SVD approximates coefficients $\inner{\cx}{u_i}$ corresponding to singular values with $\sigma_i^2\geq \al$. Coefficients corresponding to smaller singular values are set to zero. To improve the reconstruction we propose to train a CNN that reconstructs the missing part
$ \cx  - \Bo_\al \Ao \cx = \sum_{\sigma_i^2 < \al}  \inner{\cx}{u_i} u_i$ from the estimated part $\Bo_\al \cy \simeq \Bo_\al \Ao \cx$. 
Adding nonzero values to the truncated singular value components allows to better approximate non-smooth elements. 
To achieve a deep learning based expansion of singular values,
we consider a family of regularizing networks $(\Ro_\al)_{\al>0}$ of the form
\begin{equation}\label{eq:recnet}
\Ro_\al (\cy)   \triangleq
 \sum_{\sigma_i^2\geq \al} \frac{1}{\sigma_i}\inner{\cy}{v_i} u_i  +
 \nun_\al
 \kl{\sum_{\sigma_i^2\geq \al} \frac{1}{\sigma_i}\inner{\cy}{v_i} u_i}
 \; \text{ with }
\nun_\al(z)  \triangleq \sum_{\sigma_i^2 < \al}  \inner{\Uo_\al z}{u_i} u_i \,.
\end{equation}
Here $\sum_{\sigma_i^2\geq \al} {1}/{\sigma_i}\inner{\cy}{v_i} u_i
= \Bo_\al \cy$
is the truncated SVD reconstruction, $\Uo_\al  \colon \R^{N^2}  \to \R^{N^2}$ a CNN that operates on elements in $\R^{N^2}$  as $N \times N$ images, and $\al$ is the regularization parameter.

Using the projected network $\nun_\alpha$ instead of the original network $\Uo_\al$ implies that the non-vanishing coefficients of the truncated SVD expansion   are unaffected by the network. As a consequence, $\Bo_\al$ and $\Ro_\al$ reconstruct
the same low frequency part $\sum_{\sigma_i^2\geq \al} \frac{1}{\sigma_i}\inner{\cy}{v_i} u_i$.
The subsequent application of the network $\Phi_\al$ adds high frequency components that are missing in the truncated SVD reconstruction.
In particular, $\nun_\al $ is a deep CNN with $\ran(\nun_\al )\subset \operatorname{span}\{u_i \mid \sigma_i^2<\al\}$ which results in a  data driven continuation of the truncated SVD.
The CNNs $\nun_\al$ are trained to map the truncated SVD reconstruction $\Bo_\al(\cy)$ lying in the space spanned by the reliable basis elements (corresponding to large enough singular values of the operator $\Ao$) to the coefficients of the unreliable basis elements (corresponding to small singular values).
Details on the selection of the regularization parameter,
the used network architecture, and the network training
will be discussed in Section~\ref{sec:num}.

\section{Numerical  results}
\label{sec:num}

\subsection{Discretization}

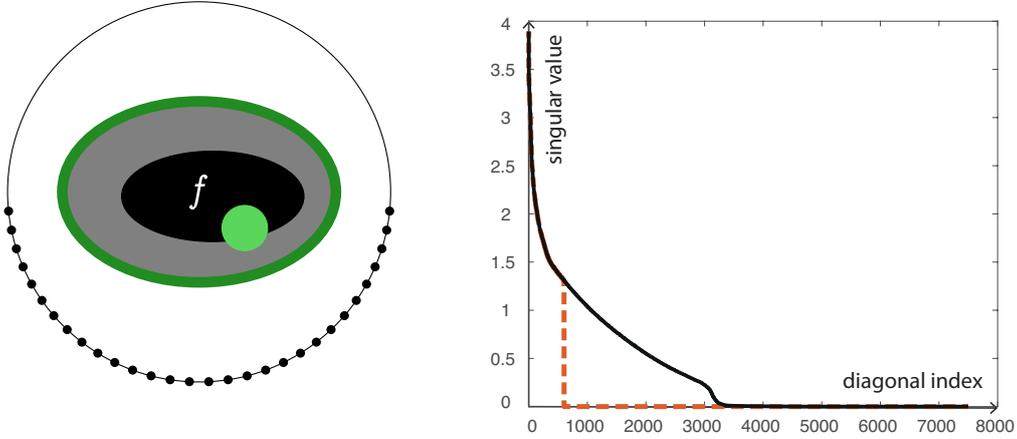
\begin{figure}[htb!]
\centering
\begin{tikzpicture}[scale=0.60]
\draw (0,6.5) node{};
\draw (0,-5.3) node{};
\filldraw[color=gray] (0,0) ellipse (3cm and 2cm);
\draw[line width=4pt,color=lime] (0,0) ellipse (3cm and 2cm);
\filldraw[color=black] (0.3,-0.10) ellipse (2cm and 1cm);
\draw[color=white] (0,0) node{{\Large $f$}};
\filldraw[color=orbl] (1,-0.8) circle(0.5);
\draw [color=black](0,0) circle(4.2);
\equic[4.2 cm]{30}{0}
\end{tikzpicture}\hspace{0.08\textwidth}
\includegraphics[width=0.48\textwidth]{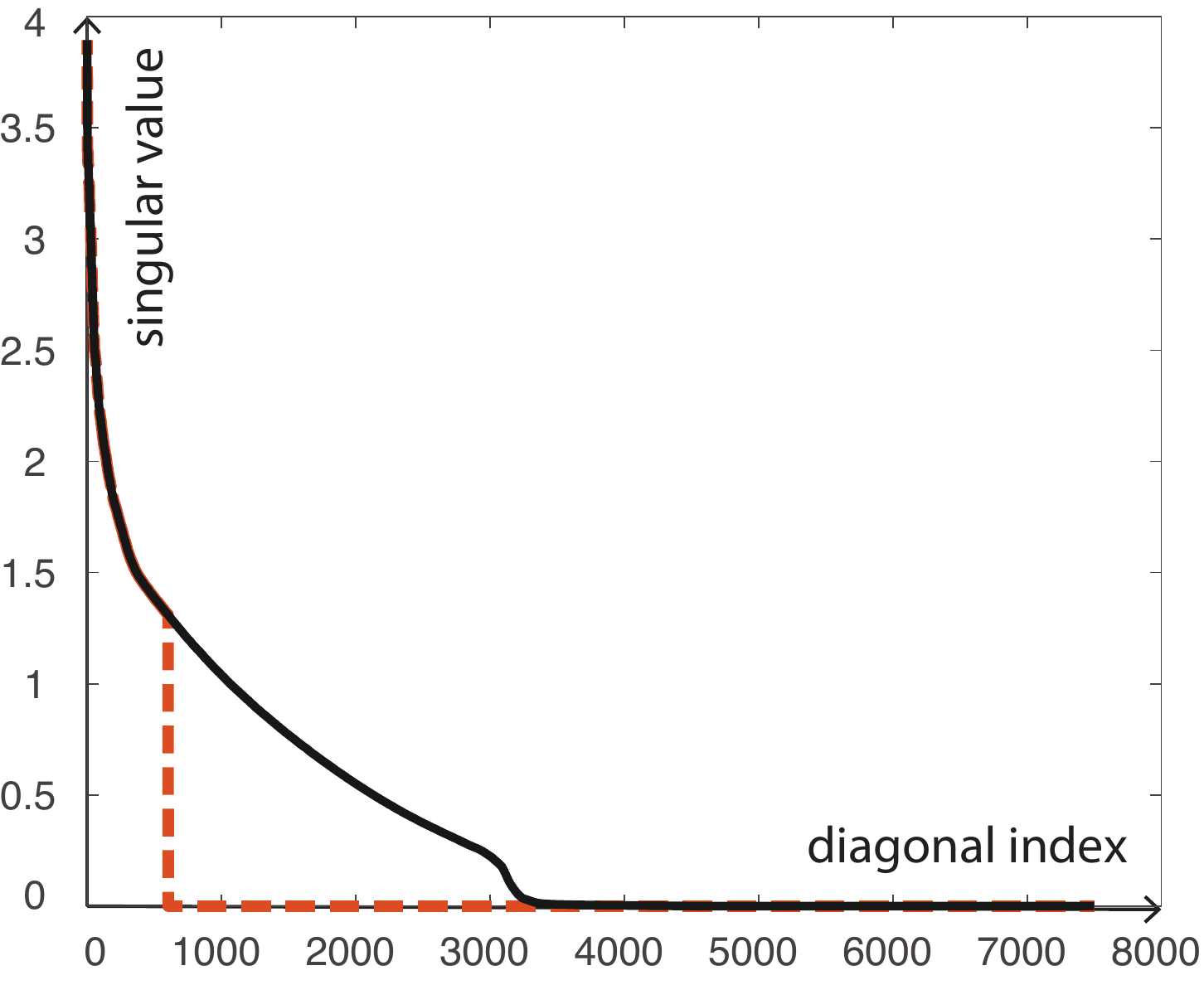}
\caption{\textbf{Considered limited view measurement geometry.}
Left: The phantom is allowed to take non-vanishing values in the square $[-1,1]$. Measurements are made on a semi circle of radius one, resulting in a severally ill-posed reconstruction
problem. Right:  solid black curve shows the singular
values of the corresponding discrete system matrix 
using $16384$ basis functions and $150400$ data points.
A large fraction of the singular values is zero or close to zero, reflecting the ill-posedness. The red dashed curve indicates the used singular values for \SI{7}{\percent} additive Gaussian noise.}\label{fig:geometry}
\end{figure}

To evaluate the regularization method described in 
Section~\ref{sec:reg}, we use a discretization of PA sources using  $128\times 128$ translated Kaiser-Bessel
basis functions $\varphi_i (\rr) = \varphi (\rr - \rr_i)$
with center points $\rr_i = (-1,1) + 2(i_1-1, i_2-1) /127 $ for
$i_1, i_2 = 1, \dots, 128$.
We have chosen support radius $a=0.055$,    window taper $\gamma=7$ and smoothness parameter $m=2$.
We consider measurement on a semi circle with radius
1, see Figure~\ref{fig:geometry}.
We use a total number of $N_s = 400$ detector positions and $N_t=376$ equidistant measurements in time in the interval $[0,3.75]$, leading to a system matrix of size $\Ao\in \R^{150400 \times 16384}$.

The entries of the system matrix $\Ao$ have been constructed by computing the solution of the initial value problem for the basis function $\varphi$ centered at $0$ for different detector-source distances. For that purpose we evaluated solution formula
\begin{equation*}
p(r,t)\triangleq\frac{1}{2\pi} \frac{\partial}{\partial t}\int_0^t\int_{\partial B_1(0)}\frac{r\varphi(\rs+r\omega)}{\sqrt{t^2-r^2}}  \ \rmd \omega \ \rmd r \,,
\end{equation*}
for varying $r = \norm{\rs}$. Using that solution of the wave equation
with initial data $\varphi_i (\rr) = \varphi (\rr-\rr_i)$ evaluated at $\rs$ equals
$p (\rs-\rr_i, t)$, the entries  $\Ao \varphi_i $ have been computed numerically by using linear interpolation. The limited view scenario results in a
severely ill-posed reconstruction problem \eqref{eq:linsys}.

\subsection{Network design}

For image reconstruction, we use the reconstruction
network \eqref{eq:recnet}, which
can be written in the form $\Ro_\alpha  = (\Id   + \Pro_{\alpha} \Uo_\al  )  \Bo_\al$.
Here $\Bo_\al$ is the truncated SVD,  $\Uo_\al \colon \R^{N^2} \to \R^{N^2} $ a trained deep neural network,
and  $\Pro_\al$  the projection onto the space
$\operatorname{span}\{ u_i \mid \sigma_i^2<\al\}$
spanned by the truncated singular vectors.
To choose the regularization parameter $\al$ we use the following  empirical approach.  We compute  noise-free data $\Ao \cx$ and noisy data $\cy = \Ao \cx + \xi$, respectively, for different PA sources and different realizations of additive Gaussian noise. Then we evaluate  the truncated SVD  reconstructions $\Bo_\al \cy$ for different regularization parameters and empirically chose such a regularization parameter, where
$\Bo_\al \Ao \cx$  is close to $ is \Bo_\al \cy$
is small. As shown in the right image in Figure~\ref{fig:geometry},
for  additive Gaussian noise with a standard deviation equals $\SI{7}{\percent}$ of the maximum
of $\Ao \cx$, this results in a regularization parameter where the first 604 singular values are kept.

To train the networks $\Ro_\al$
we generated $N_{\rm train} = 3500$ training pairs $(\cy_n, \cx_n)_{n=1}^{N_{\rm train}}$ where $\cx_n$ are training phantoms and  $\cy_n = \Ao \cx_n$ the corresponding exact data. The phantoms $\cx_n$ consist of randomly modified Shepp-Logan phantoms, where the position, size and orientation of the objects are randomly chosen and additionally some finer structures are added. Furthermore, a randomly generated deformation was applied to  the Shepp-Logan like phantoms.
Training was performed end-to-end
using noise-free data $\Ao \cx_n = \cy_n$ by
minimizing the error functional
\begin{equation*}
E_\alpha (\Uo_{\al}) =
\frac{1}{N_{\rm train}}
\sum_{n=1}^{N_{\rm train}}
\norm{ \cx_n - (\Id   + \Pro_{\alpha} \Uo_{\al}  )  \Bo_\al \Ao \cx_n  }_2^2\,,
\end{equation*}
over all free parameters in $\Uo_{\al}$.
In all cases, the network was trained for 70 epochs using stochastic gradient descent with a learning rate of 0.01 and a momentum parameter 0.99. In our results we use the
U-net architecture \cite{ronneberger2015unet} for $\Uo_{\al}$; however our
approach can also be combined with different CNN
architectures.

\begin{figure}[htb!]%
	\centering
	\subfloat{\includegraphics[trim=0 1.5cm 0 1cm,clip,width=0.49\textwidth]{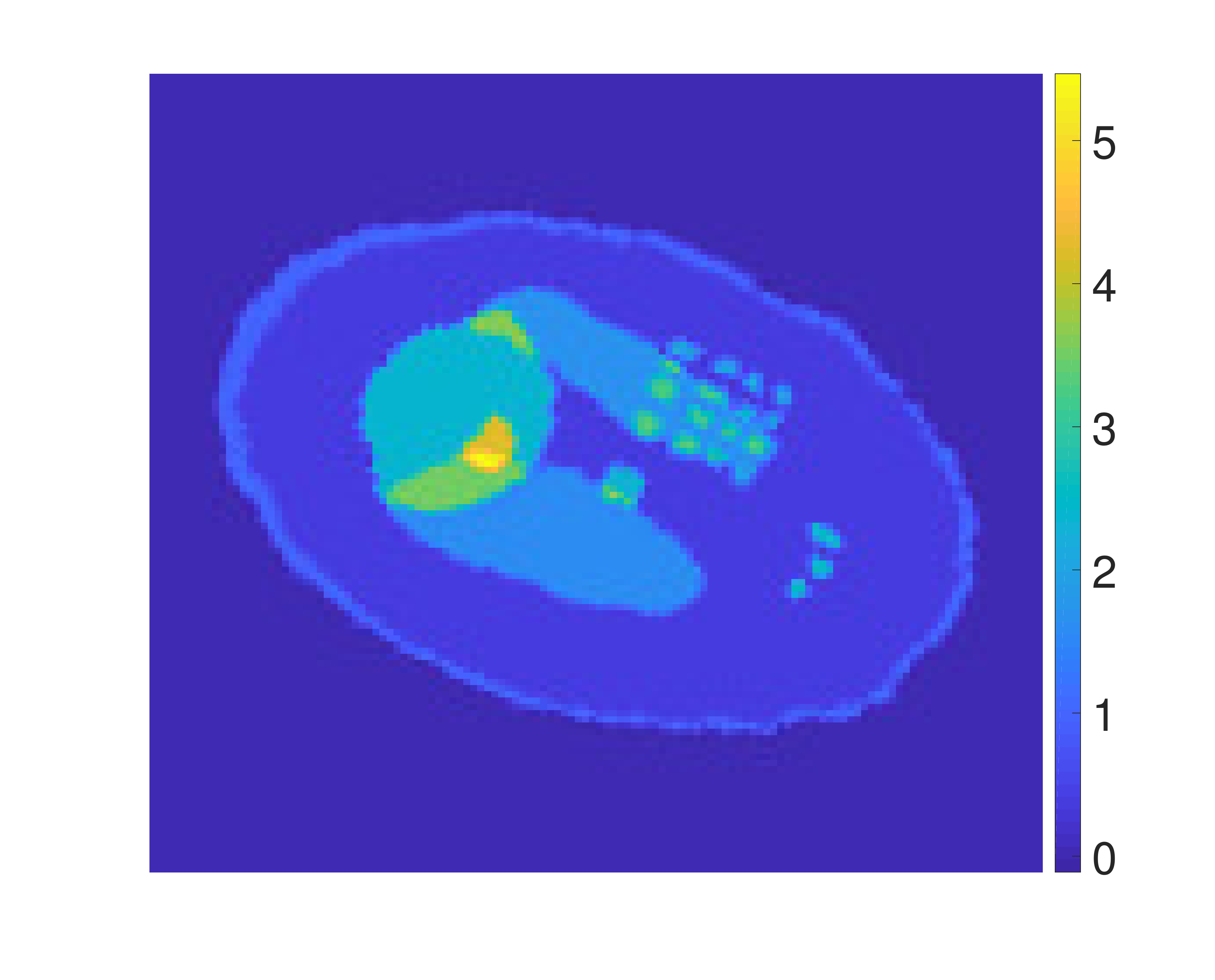}}%
	\subfloat{\includegraphics[trim=0 1.5cm 0 1cm,clip,width=0.49\textwidth]{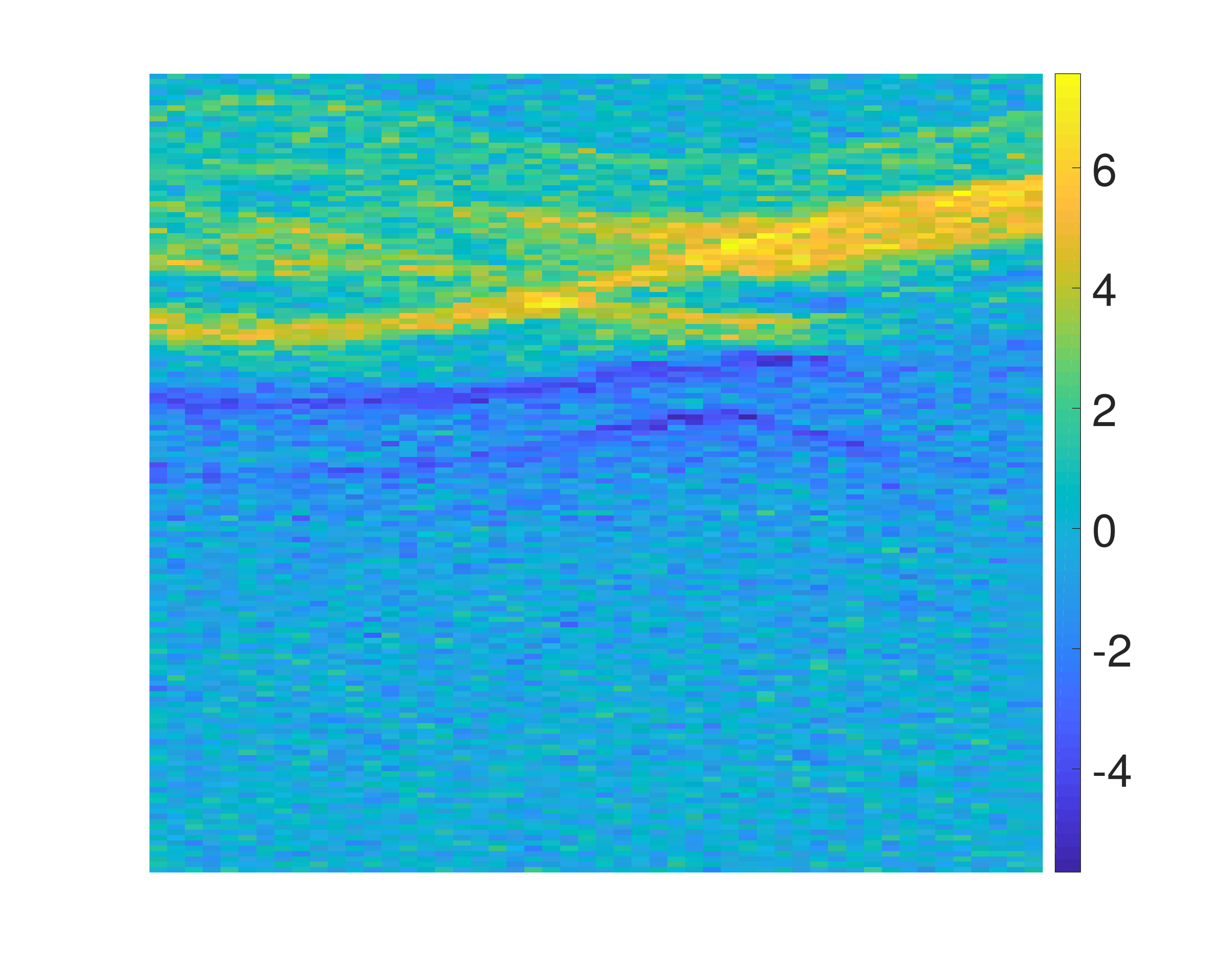}}\\
	\caption{\textbf{Phantom and noisy data.} Left: Test phantom generated by the method described above not contained in the training set. Right: Limited view data with $\SI{7}{\percent}$ additive Gaussian noise.}\label{fig:phan}
\end{figure}

\subsection{Reconstruction results}

\begin{figure}[htb!]%
	\centering
	\subfloat{\includegraphics[trim=0 1.5cm 0 1cm,clip,width=0.33\textwidth]{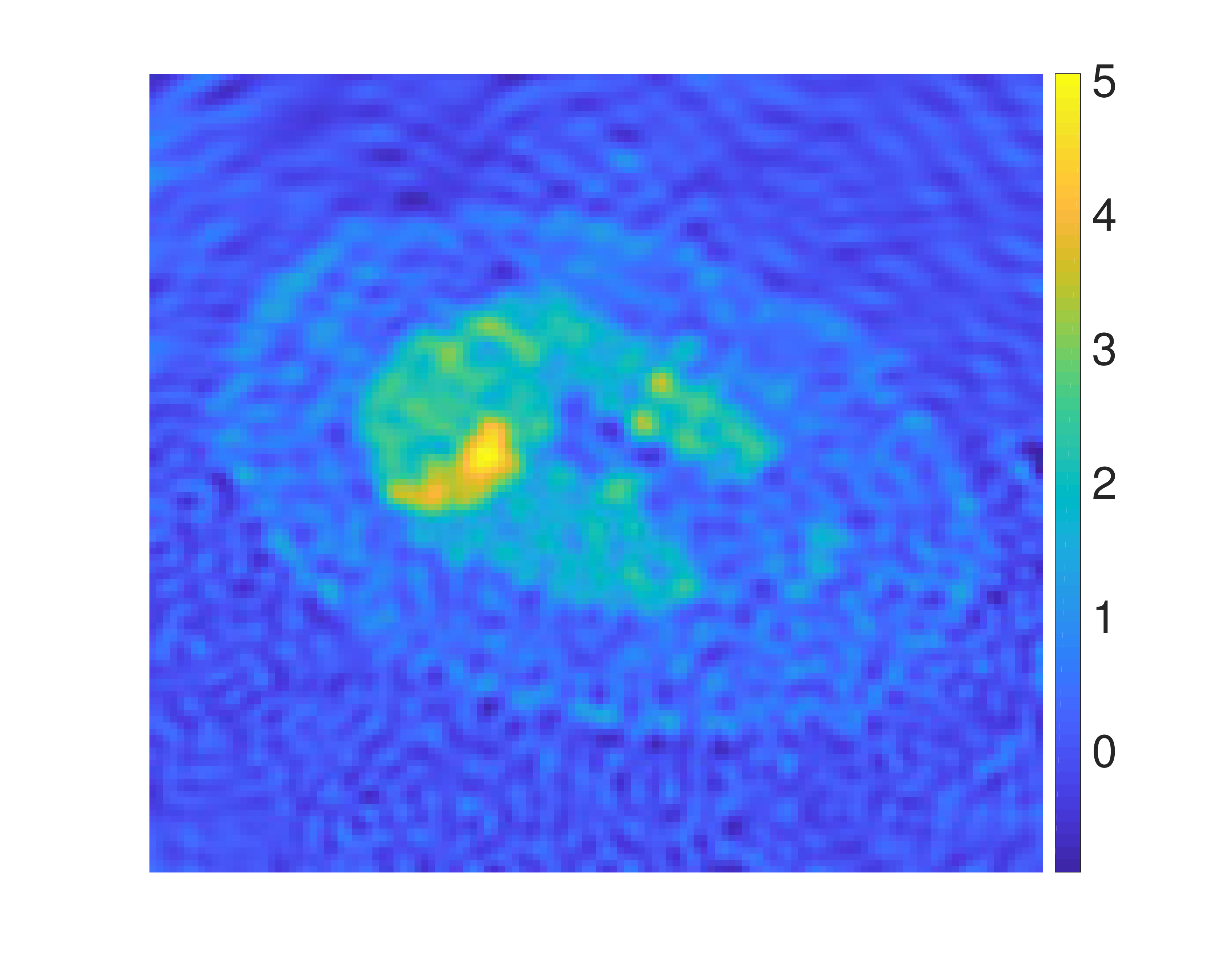}}%
	\subfloat{\includegraphics[trim=0 1.5cm 0 1cm,clip,width=0.33\textwidth]{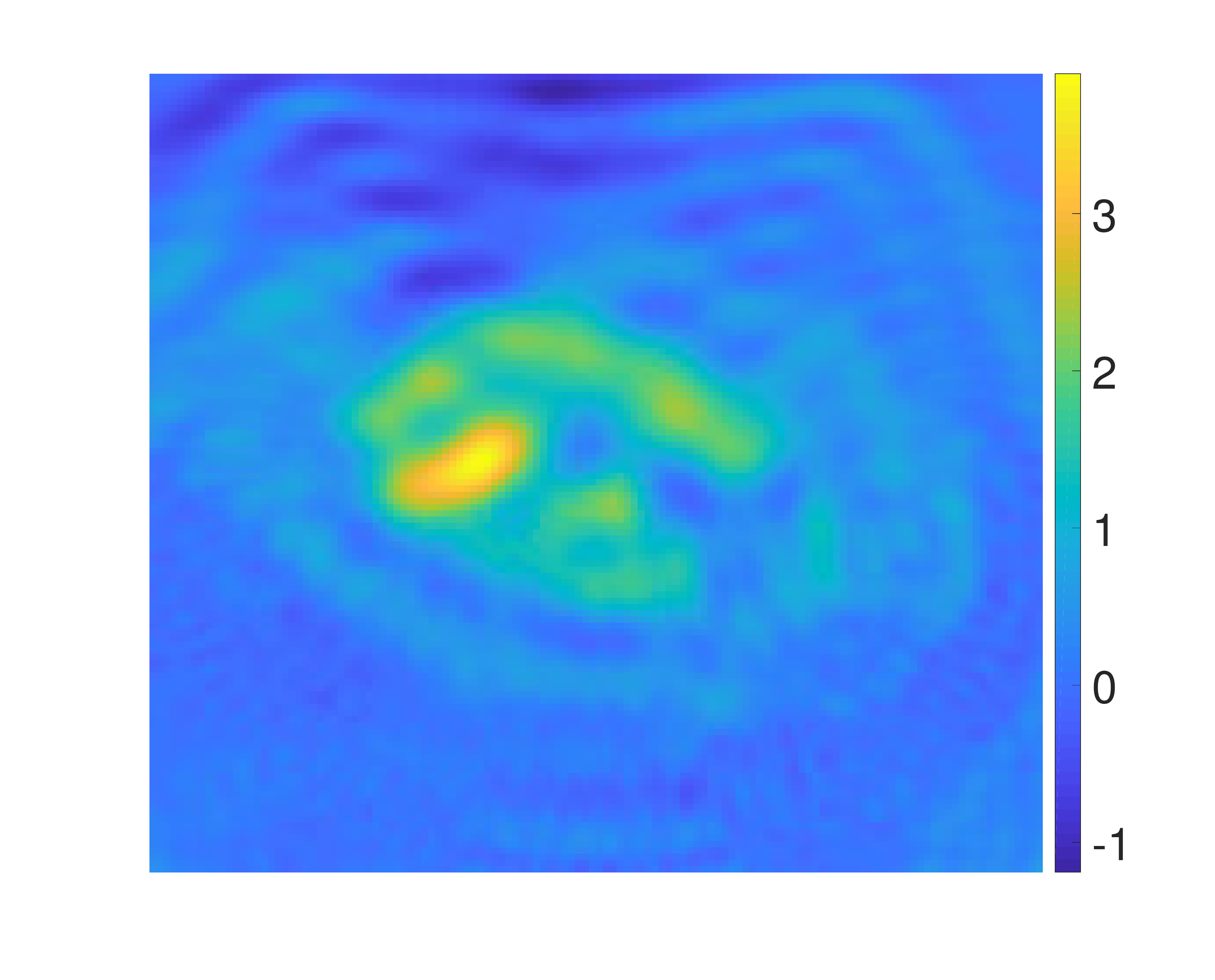}}%
	\subfloat{\includegraphics[trim=0 1.5cm 0 1cm,clip,width=0.33\textwidth]{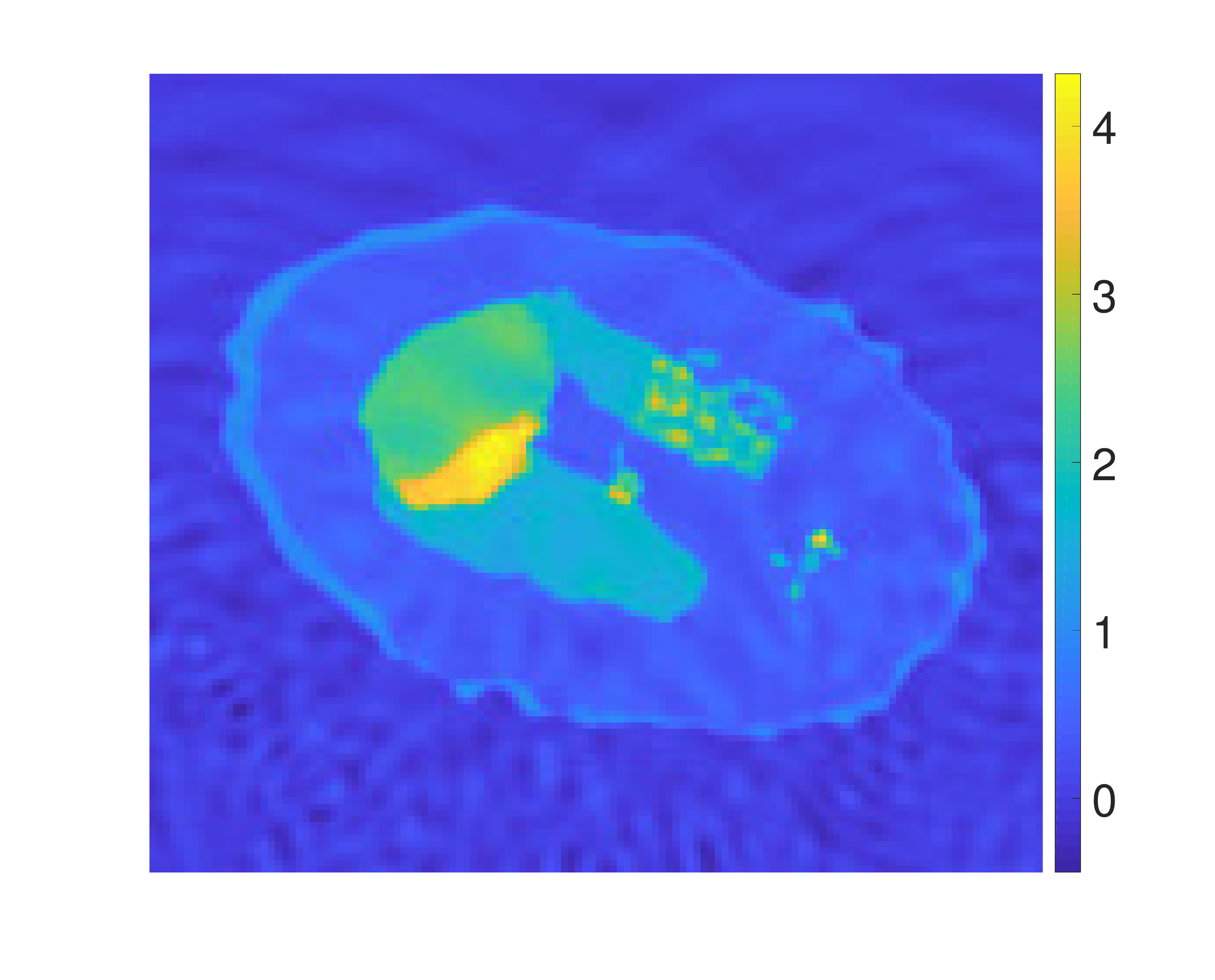}}\\[-0.5em]
	\subfloat{\includegraphics[trim=0 1.5cm 0 1cm,clip,width=0.33\textwidth]{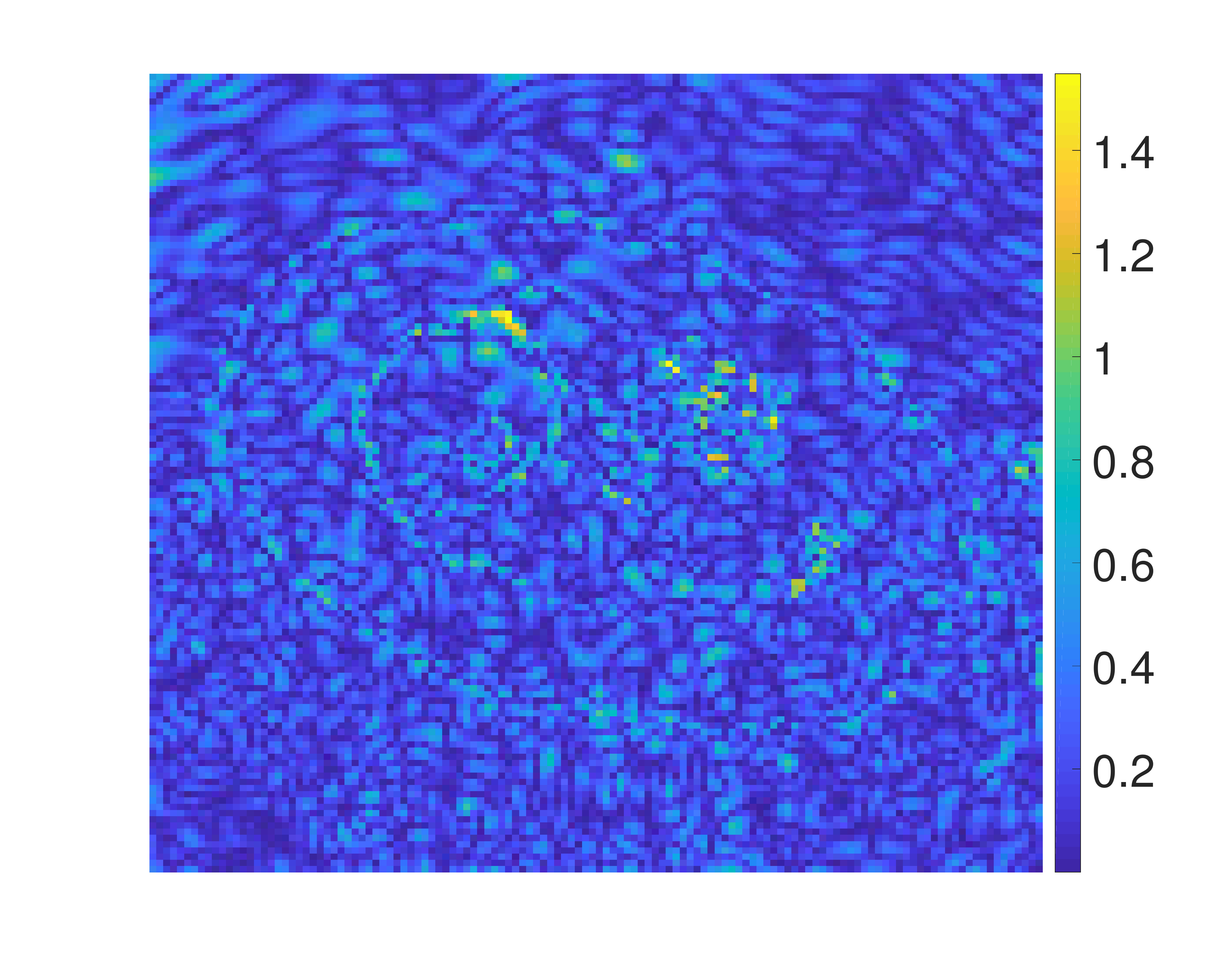}}%
	\subfloat{\includegraphics[trim=0 1.5cm 0 1cm,clip,width=0.33\textwidth]{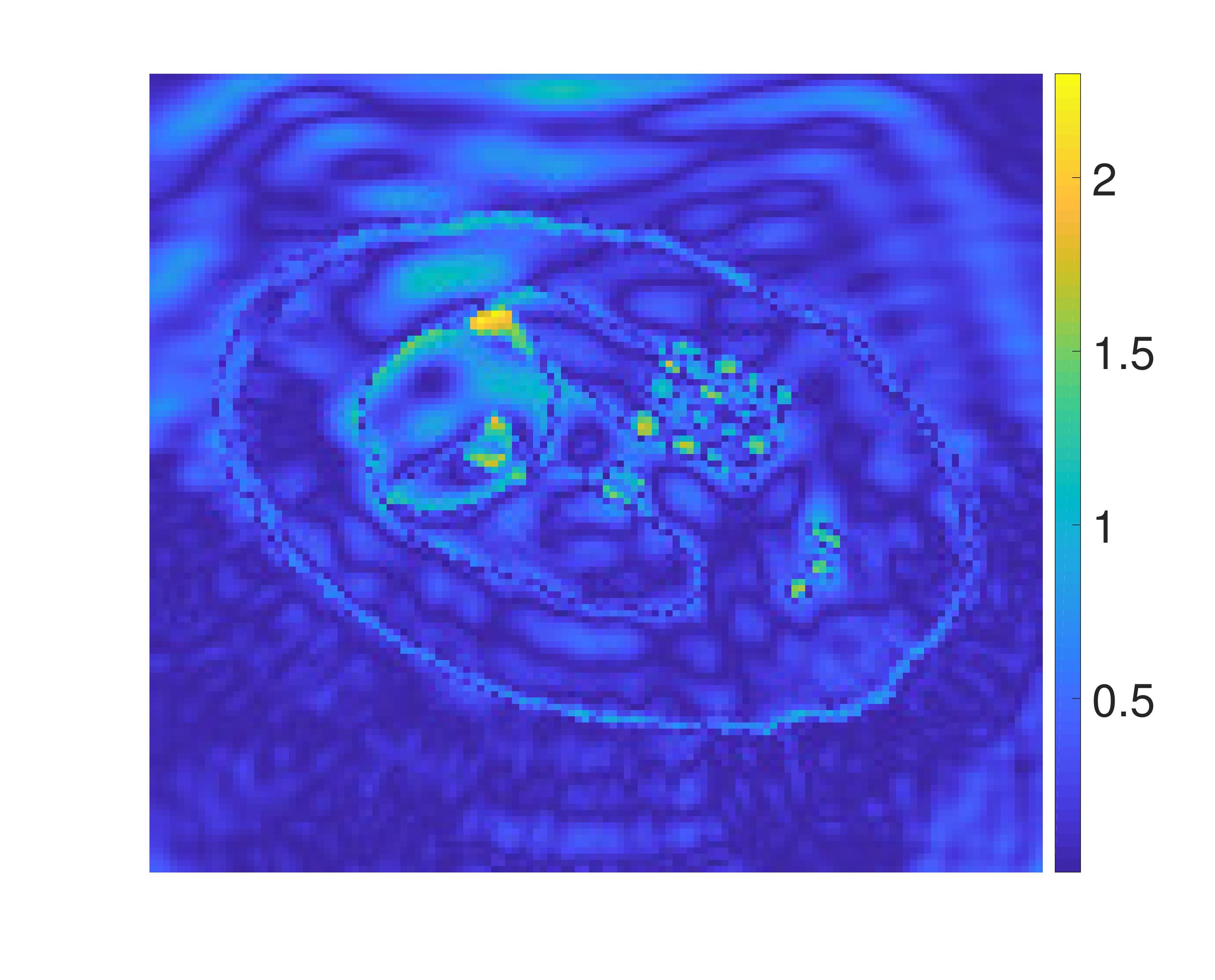}}%
	\subfloat{\includegraphics[trim=0 1.5cm 0 1cm,clip,width=0.33\textwidth]{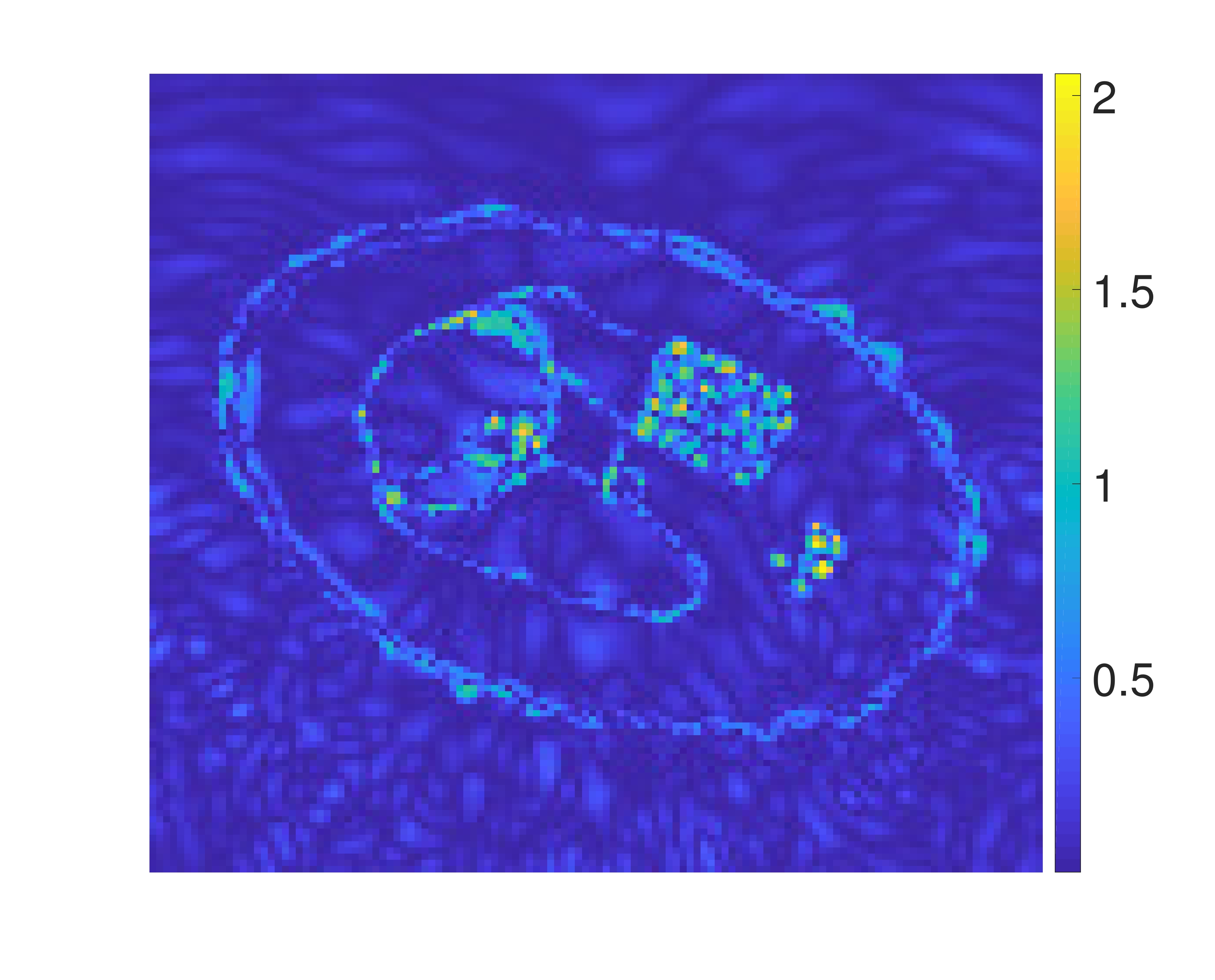}}\\
	\caption{\textbf{Reconstructions result for $\SI{7}{\percent}$ noise.} Top row: Optimally truncated SVD (left),
truncated SVD used as input to the CNN (center),
 and proposed CNN continuated SVD (right).
 Bottom row: Absolute differences from the ground truth corresponding to the reconstruction in the top row. }\label{fig:train}
\end{figure}

Figure~\ref{fig:phan} (left) shows a randomly chosen modified Shepp-Logan type phantom $\cx$. The phantom is  different from all training phantoms but generated using the same random model. The corresponding data  obtained by multiplying $\cx$ with the system matrix and adding Gaussian noise with a standard deviation of $\SI{7}{\percent}$ of the maximum of $\Ao \cx$
is shown in Figure~\ref{fig:phan}, right.
Figure~\ref{fig:train} (right column) shows the reconstruction with the proposed regularization network. The intermediate truncated
SVD reconstruction $\Bo _\al \cy$ used as input for the network $(\Id   + \Pro_{\alpha} \Uo_\al  )$  is shown in the second column. For comparison purpose, in the left column in Figure~\ref{fig:train}
we show the optimal truncated SVD reconstruction, where the regularization parameter is chosen in such a way, that $\Bo _\al \cy$  
has the minimal $\ell^2$-difference to the ground truth image.

\begin{figure}[htb!]%
	\centering
	\subfloat{\includegraphics[trim=0 1.5cm 0 1cm,clip,width=0.49\textwidth]{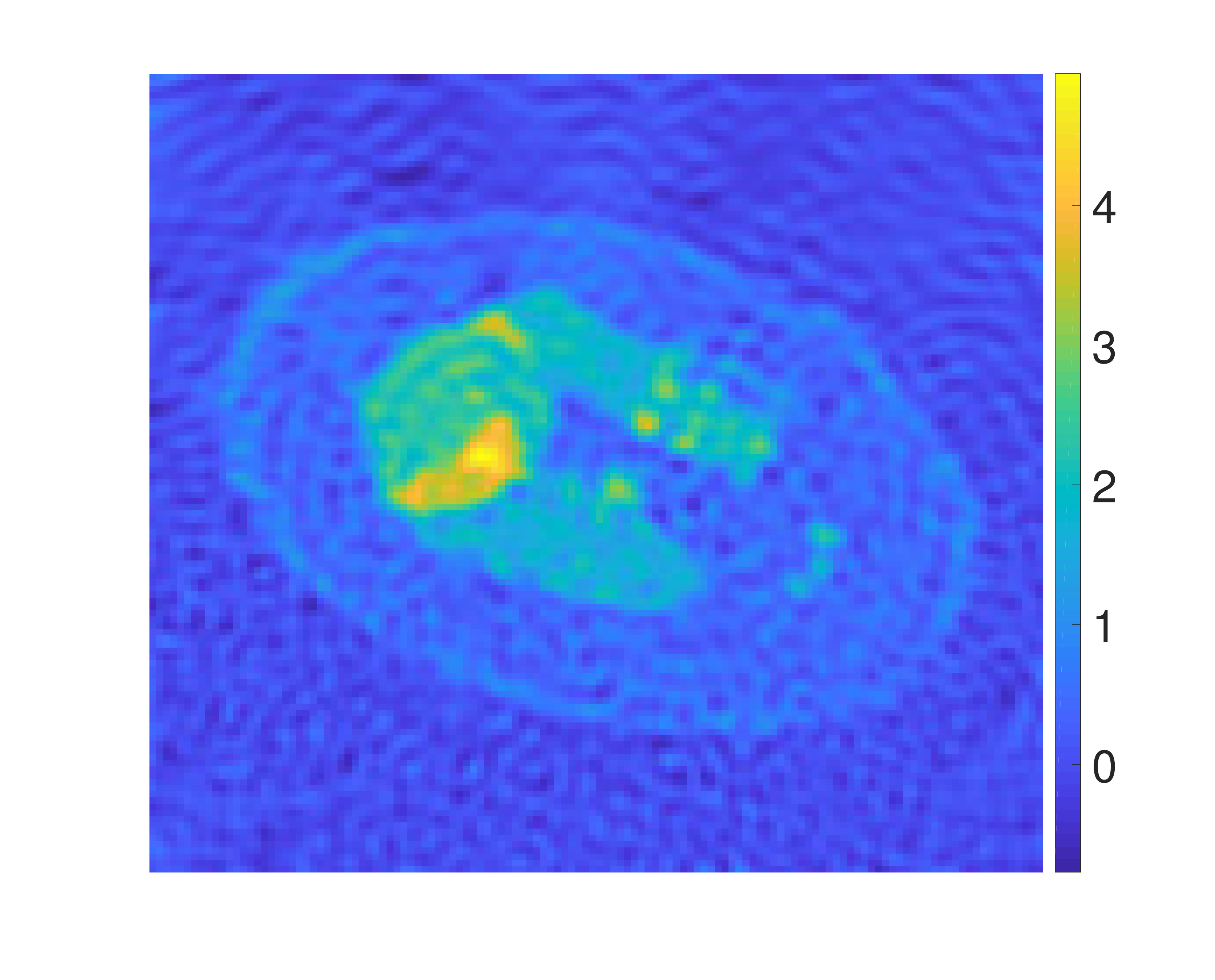}}%
	\subfloat{\includegraphics[trim=0 1.5cm 0 1cm,clip,width=0.49\textwidth]{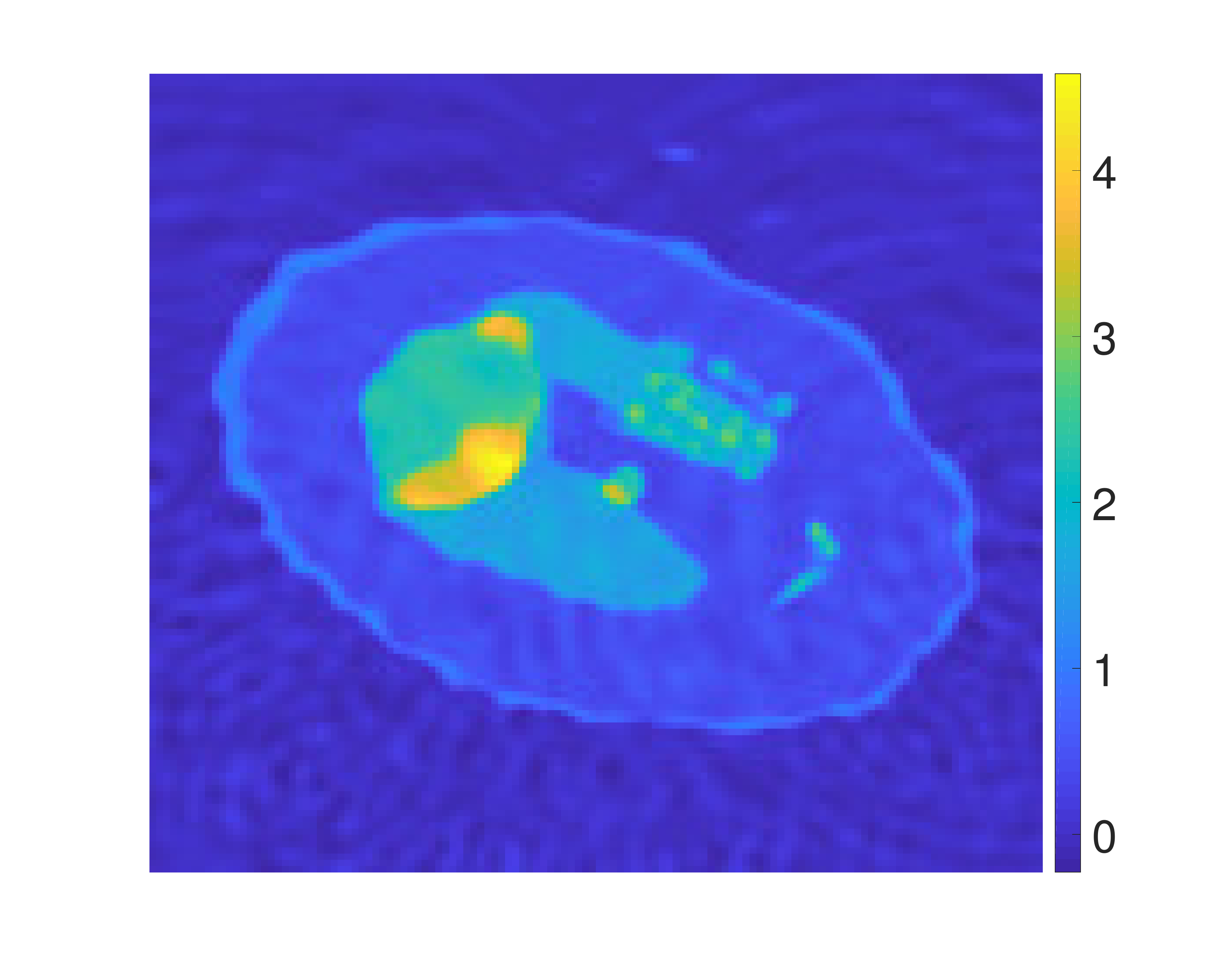}}\\[-0.5em]
	\subfloat{\includegraphics[trim=0 1.5cm 0 1cm,clip,width=0.49\textwidth]{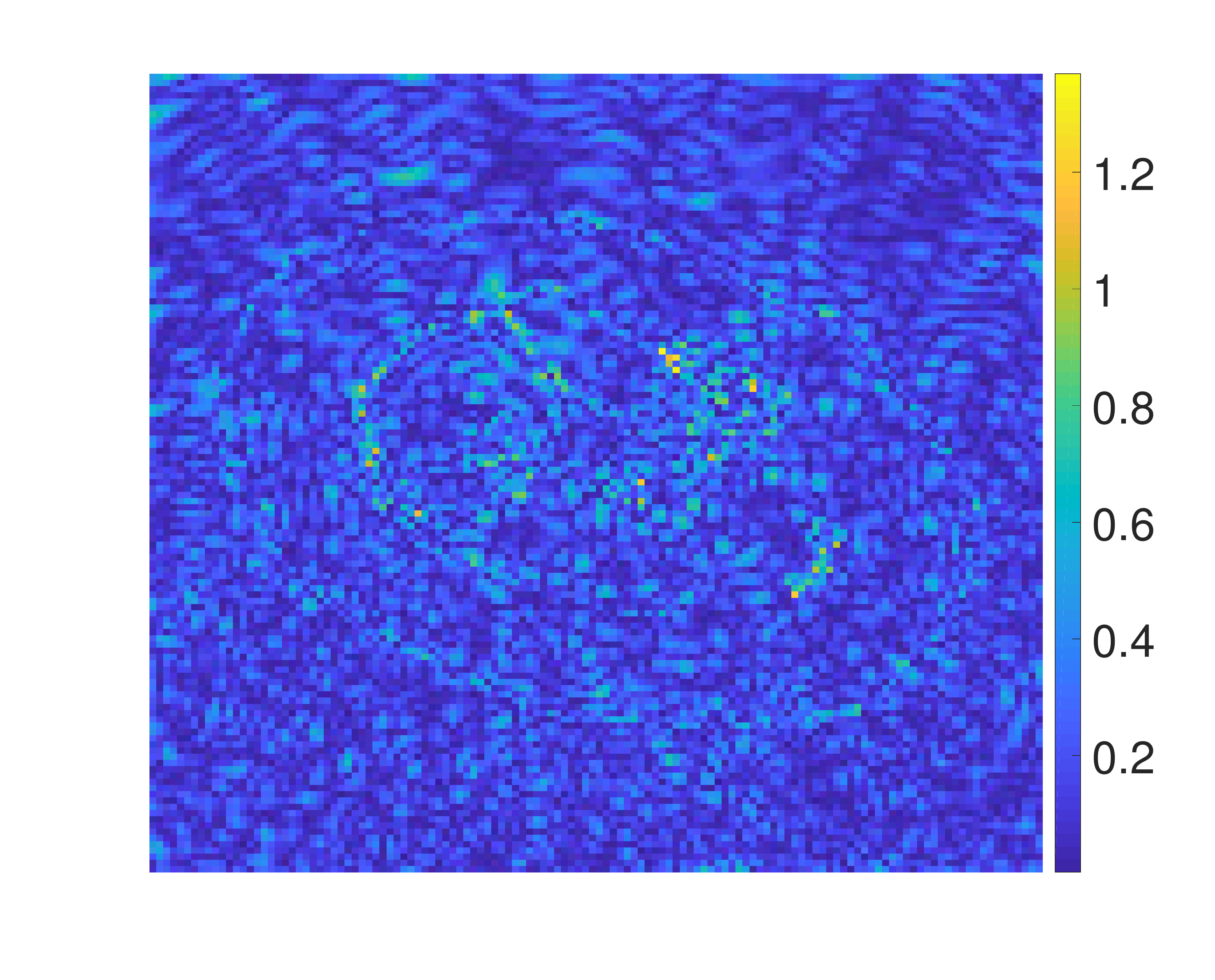}}%
	\subfloat{\includegraphics[trim=0 1.5cm 0 1cm,clip,width=0.49\textwidth]{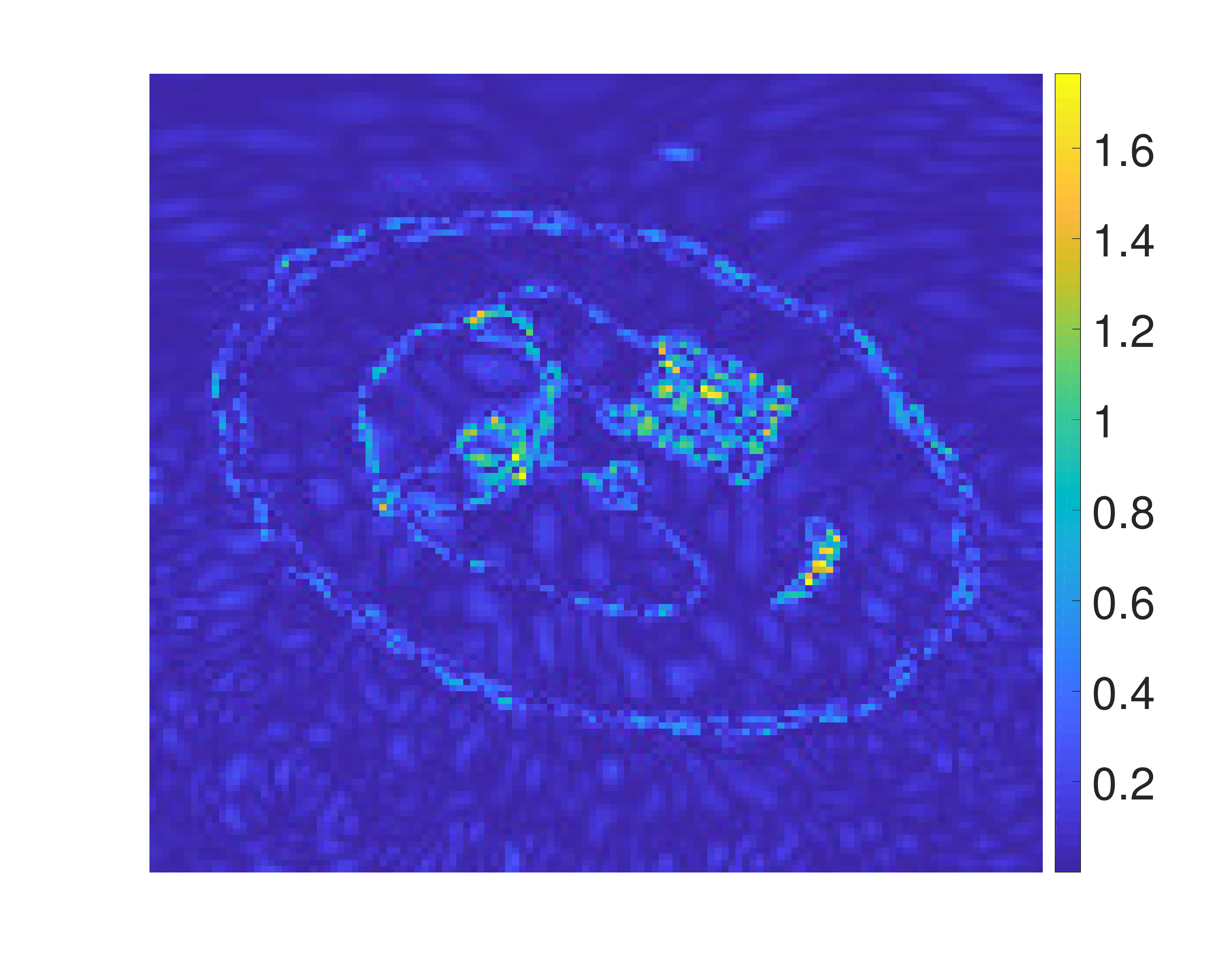}}\\
	\caption{\textbf{Reconstructions result for $\SI{4}{\percent}$ noise.} Top row: Optimally truncated SVD (left),
 and proposed CNN continuated SVD (right). Bottom row: Absolute differences from the ground truth corresponding to the reconstruction in the top row. }\label{fig:train2}
\end{figure}

We observe, that the reconstructions with the proposed regularization method are visually superior to the reconstructions using truncated SVD with optimally chosen regularization parameter. Evaluation of the relative mean squared error $
\frac{1}{500}\sum_{i=1}^{500}
{\|\Ro_\al \cy_i - \cx_i\|_2}/{\|\cx_i\|_2}
$ averaged over 500 examples with noisy data not contained in the training data set, yielded an average error of 0.0887, opposed to 0.1563 for the reconstructions coming from the optimal
truncated SVD.
Reconstruction results from data with $\SI{4}{\percent}$ additive Gaussian noise are shown in Figure~\ref{fig:train2}. Again the
proposed regularization network clearly outperforms the optimal truncated SVD reconstruction.

\section{Conclusion}
\label{sec:conclusion}

In this paper we investigates the use of regularizing networks
for the limited view problem of PAT. The truncated SVD is  used to reconstruct an intermediate
low resolution approximation, a subsequent projection
network recovers the high frequency parts.
The reconstruction networks are shown to significantly improve the truncated SVD as regularization method. Further work will be done to refine the network architecture, combine it with different regularization methods, compare it to other deep
learning based regularization methods, and to evaluate various approaches on experimental data.

\section*{Acknowledgement}

The research of  S.A. and M.H. has been  supported by the Austrian Science Fund (FWF), project P 30747-N32; the work of R.N. has been supported by the FWF, project P 28032. We acknowledge the support of NVIDIA Corporation with the donation of the Titan Xp GPU used for this research.

\end{document}